\documentclass[preprint,12pt]{elsarticle}

\usepackage{amssymb, url}
\usepackage{lscape}
\usepackage{color}
\usepackage{graphicx}

\journal{}

\begin{document}

\begin{frontmatter}



\title{Quasi-Monte Carlo point sets with small $t$-values and WAFOM} 


\author[harase]{Shin Harase\corref{cor1}}
\ead[harase]{harase@craft.titech.ac.jp}
\address[harase]{Graduate School of Innovation Management, Tokyo Institute of Technology, 
W9-115, 2-12-1 Ookayama, Meguro-ku, Tokyo, 152-8550, Japan.}

\cortext[cor1]{Corresponding author (Tel: +81-3-5734-3517)}

\begin{abstract}
The $t$-value of a $(t, m, s)$-net is an important criterion of point sets for quasi-Monte Carlo integration, 
and many point sets are constructed in terms of \textcolor{blue}{the} $t$-values, as 
this lead\textcolor{blue}{s} to small integration error bounds.
Recently, Matsumoto, Saito, and Matoba proposed the Walsh figure of merit (WAFOM) as a quickly computable criterion of point sets that ensure\textcolor{blue}{s} higher order
 convergence for function classes of very high smoothness. 
In this paper, we consider a search algorithm for point sets whose $t$-value and WAFOM are both small\textcolor{blue}{,} 
so as to be effective for a wider range of function classes. 
For this, we fix digital $(t, m, s)$-nets with small $t$-values (e.g., Sobol' or Niederreiter--Xing nets) in advance, 
apply random linear scrambling, and select scrambled digital $(t, m, s)$-nets in terms of WAFOM. 
Experiments show that the \textcolor{blue}{resulting} point sets improve the rates of convergence for smooth functions 
and are robust for non-smooth functions. 
\end{abstract}


\end{frontmatter}


\newtheorem{theorem}{Theorem}
\newtheorem{lemma}[theorem]{Lemma}
\newtheorem{corollary}[theorem]{Corollary}
\newdefinition{definition}{Definition}
\newdefinition{remark}{Remark}
\newproof{proof}{Proof}
\newtheorem{proposition}[theorem]{Proposition}


\section{Introduction}
For a Riemann integrable function $f:[0,1)^s \to \mathbb{R}$,
we consider the integral $\int_{[0,1)^s} f(\mathbf{x}) \textrm{d} \mathbf{x}$ 
and its approximation by quasi-Monte Carlo integration: 
\begin{eqnarray} \label{eqn:monte carlo}
 \int_{[0,1)^s} f(\mathbf{x}) \textrm{d} \mathbf{x} \approx \frac{1}{N}\sum_{k = 0}^{N-1} f(\mathbf{x}_k),
\end{eqnarray}
where the point set $P :=\{ \mathbf{x}_0, \ldots, \mathbf{x}_{N-1} \} \subset [0,1)^s$ is chosen deterministically.

A typical quasi-Monte Carlo point set $P$ is a low-discrepancy point set based on the $t$-value of a $(t, m, s)$-net. 
Thus, the $t$-value is probably the most important criterion of quasi-Monte Carlo point sets \cite{MR3038697, MR2683394, MR1172997}. 

Matsumoto, Saito, and Matoba \cite{MSM} recently proposed 
the Walsh figure of merit (WAFOM) as another criterion of 
quasi-Monte Carlo point sets \textcolor{blue}{to ensure} higher order convergence 
for function classes of very high smoothness.
WAFOM is also quickly computable, 
and this efficiency enables us to search for quasi-Monte Carlo point sets using a random search. 
\textcolor{blue}{From an analogy to} coding theory, 
since a random search is easier than a mathematical construction (e.g., the success of low-density parity-check codes), 
Matsumoto et al.\ also searched for point sets at random by minimizing WAFOM. 
In the same spirit, Harase and Ohori \cite{HO} searched for low-WAFOM point sets with extensibility (i.e., the number of points may be increased while the existing points are retained). 
In numerical experiments, these point sets are significantly effective for low-dimensional smooth functions. 
In fact, as shown later (in Remark~\ref{Remark: t-value}), low-WAFOM point sets based on a simple random search
do not always have small $t$-values in the framework of $(t, m, s)$-nets, and such point sets are sometimes inferior to classical $(t, m, s)$-nets for non-smooth functions. 

In this paper, we search for point sets whose $t$-value and WAFOM are both small, 
so as to be effective for a wider range of function classes, i.e., 
point sets combining the advantages of good $(t, m, s)$-nets and low-WAFOM point sets.
For this, we fix suitable digital $(t, m, s)$-nets (e.g., Sobol' or Niederreiter--Xing nets) in advance
and apply random linear scrambling with non-singular lower triangular matrices 
that preserves \textcolor{blue}{the} $t$-values. 
\textcolor{blue}{The key to our approach} is to select good point sets from the scrambled digital $(t, m, s)$-nets in terms of WAFOM. 
Our numerical experiments show that the obtained point sets improve 
the rates of convergence for smooth functions and are robust for non-smooth functions. 

The rest of this paper is organized as follows. 
In Section~2, we briefly recall
the definitions of digital $(t, m, s)$-nets and WAFOM. 
Section~3 is devoted to our main result: a search for low-WAFOM point sets with small $t$-values using linear scrambling.
In Section~4, we \textcolor{blue}{compare our new point sets with} other quasi-Monte Carlo point sets 
by using the Genz test function package \cite{Genz1984,Genz1987}.
Section~5 concludes the paper with some directions for future research. 

\section{Notations}
\subsection{Digital $(t, m, s)$-nets}
\textcolor{blue}{We briefly recall the definition of digital $(t, m, s)$-nets.  
Throughout this paper, we consider only the digital $(t, m, s)$-nets in base $2$.} 
Let $s$ and $n$ be positive integers. Let $\mathbb{F}_2:=\{ 0,1\}$ be the two-element field, 
and $V := \mathbb{F}_2^{s \times n}$ the set of $s \times n$ matrices. 
Let us denote $\mathbf{x} \in V$ by $\mathbf{x}:=(x_{i, j})_{1 \leq i \leq s, 1 \leq j \leq n}$ with $x_{i, j} \in \mathbb{F}_2$.  
We identify $\mathbf{x} \in V$ with the $s$-dimensional point 
\[ (\sum_{j=1}^{n} x_{1, j}2^{-j}+2^{-n-1}, \ldots, \sum_{j=1}^{n} x_{s, j}2^{-j}+2^{-n-1}) \in [0,1)^s.\]
Note that $n$ corresponds to the precision. Note also that the points are shifted by $2^{-n-1}$ 
because we will later consider WAFOM (see \cite[Remark 2.2]{MSM}). 
To construct $P:=\{ \mathbf{x}_0, \mathbf{x}_1,\ldots, \mathbf{x}_{2^m-1} \} \subset [0,1)^s$, 
we often use the following construction scheme called the {\it digital net}. 

\textcolor{blue}{
\begin{definition}[Digital net]
Consider $n \times m$ matrices $C_1, \ldots, C_s \in \mathbb{F}_2^{n \times m}$. 
For $h = 0, 1, \ldots, 2^m-1$, let $h = \sum_{l = 0}^{m-1} h_l 2^l$ with $h_l \in \mathbb{F}_2$ 
be the expansion of $h$ in base $2$. We set $\mathbf{h} := {}^t(h_0, \ldots, h_{m-1}) \in \mathbb{F}_2^m$, where ${}^t$ represents the transpose. 
We set $\mathbf{x}_h := {}^t(C_1 \mathbf{h},\ldots, C_s \mathbf{h}) \in V$. 
Then, the point set $P:= \{ \mathbf{x}_0, \ldots, \mathbf{x}_{2^m-1} \}$ is called a {\it digital net} over $\mathbb{F}_2$ 
and  $C_1, \ldots, C_s$ are the {\it generating matrices} of the digital net $P$. 
\end{definition}}
Throughout this paper, we assume $P$ is a digital net. 
Note that $P \subset V$ is an $\mathbb{F}_2$-linear subspace of $V$. 

\begin{definition}[$(t, m, s)$-net]
Let $s \geq 1$, \textcolor{blue}{and let} $0 \leq t \leq m$ be integers. 
Then, a point set $P$ consisting of $2^m$ points in $[0, 1)^s$ is called a {\it $(t, m, s)$-net} (in base $2$) 
if every subinterval $J = \prod_{i = 1}^s[a_i 2^{-d_i}, (a_i +1) 2^{-d_i})$ in $[0, 1)^s$ with integers $d_i \geq 0$ and $0 \leq a_i < 2^{d_i}$ 
for $1 \leq i \leq s$ and of volume $2^{t-m}$ contains exactly $2^t$ points of $P$. 
\end{definition}

\textcolor{blue}{
\begin{definition}[$t$-value]
If $t$ is the smallest value such that $P$ is a $(t,m, s)$-net, then we call this the $t$-value (or exact quality parameter). 
\end{definition}
\begin{definition}[Digital $(t, m, s)$-net]
If $P$ is a digital net and a $(t, m, s)$-net, it is called a {\it digital $(t, m, s)$-net}.
\end{definition}
}


As a criterion, $P$ is well distributed if the $t$-value is small. 
In this framework, from the Koksma--Hlawka inequality and estimation of star-discrepancies,
the upper bound on the absolute error of (\ref{eqn:monte carlo}) is $O(2^t(\log N)^{s-1}/N)$ 
(see \cite{MR2683394, MR1172997} for details). 
There are many studies on \textcolor{blue}{the} generating matrices of digital $(t, m, s)$-nets, e.g.,  
Sobol' nets \cite{MR0219238}, Niederreiter nets \cite{MR1172997}, and Niederreiter--Xing nets \cite{MR1358190}. 
There are also some algorithms for computing \textcolor{blue}{the $t$-value} of digital nets \cite{MR3085113, MR1881672}.


\subsection{WAFOM}

Matsumoto et al.\ \cite{MSM} proposed 
WAFOM as a computable criterion of 
quasi-Monte Carlo point sets constructed by digital nets $P$. 
WAFOM has the \textcolor{blue}{potential} to ensure higher order convergence than $O(N^{-1})$ for function classes of very high smoothness (so-called {\it $n$-smooth functions}). 
In a recent talk, 
Yoshiki \cite{Yoshiki2014} modified the definition of WAFOM resulting in a more explicit upper bound for integration errors \textcolor{blue}{(see also Section~7 of \cite{MO2014})}. 
Thus, throughout this paper, 
we adopt \textcolor{blue}{his new result as our WAFOM value with some abuse of notation}. 

\textcolor{blue}{
\begin{definition}[WAFOM]
Let $P \subset V$ be a digital net. For $A = (a_{i, j}), B = (b_{i, j})\in V$, we define the inner product as 
$\langle A ,B \rangle := \sum_{1 \leq i \leq s, 1 \leq j \leq n} a_{i, j} b_{i,j} \in \mathbb{F}_2$.
For an $\mathbb{F}_2$-linear subspace $P$, let us define its perpendicular space by 
$P^{\perp} := \{ A \in V \ | \ \langle B, A \rangle = 0 \mbox{ for all } B \in P \}$.
The WAFOM (Walsh figure of Merit) of $P$ is defined by 
\[ \mbox{WAFOM}(P) := \sum_{A \in P^{\perp} \backslash \{ \mathbf{0} \} } 2^{-\mu'(A)}, \]
where we set the weight
\begin{eqnarray} \label{eqn:Yoshiki_weight}
 \mu'(A) := \sum_{1 \leq i \leq s, 1 \leq j \leq n} (j + 1) \times  a_{i,j} \quad \mbox{ for } A = (a_{i,j}) \in P^{\perp}.
\end{eqnarray}
\end{definition}
In the original definition of WAFOM, Matsumoto et al.\ \cite{MSM} considered the weight $\mu(A) := \sum_{1 \leq i \leq s, 1 \leq j \leq n} j \times a_{i, j}$ 
instead of (\ref{eqn:Yoshiki_weight}). (The weight $\mu$ was originally proposed by Dick \cite{MR2346374, MR2391005} and is now called the Dick weight.) 
Further, by replacing $c(A) := 2^{-\mu(A)}$ by $c(A) := 2^{-\mu'(A)}$ in Theorem~4.1 and Corollary~4.2 of \cite{MSM} and their proofs,
we obtain the following efficiently computable formula:
\begin{eqnarray} \label{eqn:WAFOM}
\mbox{WAFOM}(P) = \frac{1}{|P|} \sum_{\mathbf{x} \in P} \left\{ \prod_{1 \leq i \leq s} \prod_{1 \leq j \leq n} (1+(-1)^{x_{i, j}} 2^{-(j+1)}) -1 \right\}. 
\end{eqnarray}
Thus, this criterion is computable in $O(nsN)$ arithmetic operations, where $N := |P|$, 
and is computable in $O(sN)$ steps when using look-up tables (see \cite{HO}). 
}

Next, we recall the $n$-digit discretization $f_n$ of $f$ by following \cite[Section~2]{MSM}. 
For $\mathbf{x}=(x_{i, j})_{1 \leq i \leq s, 1 \leq j \leq n} \in V$, 
we define the $s$-dimensional subinterval $\mathbf{I}_{\mathbf{x}} \subset [0,1)^S$ by 
\[ \mathbf{I}_\mathbf{x} := [\sum_{j=1}^{n} x_{1, j}2^{-j}, \sum_{j=1}^{n} x_{1, j}2^{-j} + 2^{-n}) \times \cdots \times
[\sum_{j=1}^{n} x_{s, j}2^{-j}, \sum_{j=1}^{n} x_{s, j}2^{-j} + 2^{-n}).\] 
For a Riemann integrable function $f: [0,1)^s \to \mathbb{R}$, we define its $n$-digit discretization 
$f_n: V \to \mathbb{R}$ by $f_n(\mathbf{x}) := (1/{\rm Vol}(\mathbf{I}_{\mathbf{x}})) 
\int_{\mathbf{I}_{\mathbf{x}}} f(\mathbf{x}) \textrm{d} \mathbf{x}$.
This is the average value of $f$ over $\mathbf{I}_{\mathbf{x}}$. 
When $f$ is Lipschitz continuous, 
it can be shown \cite{MSM} that the discretization error between $f$ and $f_n$ on $\mathbf{I}_{\mathbf{x}}$ is 
negligible if $n$ is sufficiently large (e.g., when $n \geq 30$). 
Thus, for such $f: [0, 1)^s \to \mathbb{R}$ and large $n$, 
we may consider $({1}/|P|) \sum_{\mathbf{x} \in P} f(\mathbf{x}) \approx (1/{|P|}) \sum_{\mathbf{x} \in P} f_n(\mathbf{x})$.

Here, we assume that $f$ is an $n$-smooth function (see \cite{MR2391005} and \cite[Ch.~14.6]{MR2683394} for the definition).  
Yoshiki \cite{Yoshiki2014} gave the following Koksma--Hlawka type inequality by improving Dick's inequality (\cite[Section~4.1]{MR2743889} and \cite[(3.7)]{MSM}): 
\begin{eqnarray} \label{ineq:yoshiki}
\left| \int_{[0,1)^s} f(\mathbf{x}) \textrm{d} \mathbf{x}-\frac{1}{|P|} \sum_{\mathbf{x} \in P} f_n(\mathbf{x}) \right| \leq  \sup_{0 \leq N_1, \ldots, N_s \leq n} 
|| f^{(N_1, \ldots, N_s)}||_{\infty} \cdot \mbox{WAFOM}(P),
 \end{eqnarray}
where $||f||_{\infty}$ is the infinity norm of $f$ 
and $f^{(N_1, \ldots, N_s)} := \partial^{N_1 + \cdots + N_s}f / \partial x_1^{N_1} \cdots \partial x_s^{N_s}$. 

\textcolor{blue}{
\begin{remark}
More precisely, Yoshiki \cite{Yoshiki2014} proved an upper bound on the Walsh coefficient of wavenumber $\mathbf{k} :=(k_1, \ldots, k_s)$ as follows: 
\[ | \hat{f} (\mathbf{k}) | \leq 2^{-\mu' (\mathbf{k})} || f^{(N_1, \ldots, N_s)} ||_{\infty}, \]
where $f$ is an $n$-smooth function and $k_i = \sum_{j = 1}^{N_i} 2^{a_{i, j}}$ such that $a_{i, 1} > \ldots > a_{i, N_i}$ for each $j$. 
From a similar argument as that for the proof of Theorem~3.4 and formula (3.5) in \cite{MSM}, the discretized upper bound (\ref{ineq:yoshiki}) is obtained.
\end{remark}
}

\begin{remark}
Following the discussions in \cite{MR3145585, Ohori2015, Suzuki2014, Yoshiki}, 
the best (i.e., smallest) value of $\log({\rm WAFOM}(P))$ is $O(-m^2/s)$ for $P$ with $|P| = 2^m$.
Thus, WAFOM can be used to search \textcolor{blue}{for} a digital net $P$ 
with higher order convergence than $O(N^{-1})$ for $n$-smooth functions. 
\end{remark}

\section{Scrambling methods} \label{Sec:scrambling}

In previous works, Matsumoto et al.\ \cite{MSM} and Harase and Ohori \cite{HO} 
searched for low-WAFOM point sets using only WAFOM as a criterion. 
In fact, the point sets obtained in these ways do not always have small $t$-values as $(t, m, s)$-nets. 
In this section, we take into account the \textcolor{blue}{$t$-value}, and search for low-WAFOM point 
sets with small $t$-values. For this, we consider the following transformation, known as {\it linear scrambling} \textcolor{blue}{\cite{MR1659004}, 
which is a subclass of (non-linear) {\it scrambling} with general permutations proposed by Owen \cite{MR1445791}.} 
\begin{proposition}
Let $C_1, \ldots, C_s \in \mathbb{F}_2^{n \times m}$ be generating matrices of a digital $(t, m, s)$-net. 
Let $L_1, \ldots, L_s \in \mathbb{F}_2^{n \times n}$ be non-singular lower triangular matrices. 
Then, the digital net with generating matrices $L_1C_1, \ldots, L_sC_s \in \mathbb{F}_2^{n \times m}$ is also a $(t, m, s)$-net. 
\end{proposition}
\textcolor{blue}{The proof is easily obtained from Theorem~4.28 in \cite{MR1172997} or Theorem~4.52 in \cite{MR2683394}.} 
Linear scrambling preserves the \textcolor{blue}{$t$-value}, so  
we cannot distinguish whether \textcolor{blue}{the} scrambled nets are good using \textcolor{blue}{the $t$-value itself}. 
Here, WAFOM can be applied to \textcolor{blue}{assess the} linearly scrambled digital $(t, m, s)$-nets.
Our algorithm proceeds as follows: 
\begin{enumerate}
\item Fix a digital $(t, m, s)$-net with a small $t$-value in advance.
\item Generate $L_1, \ldots, L_s$ at random $M$ times, and construct $P$ from $L_1C_1, \ldots, \textcolor{blue}{L_s C_s}$.
\item Select the point set $P$ with the smallest $\mbox{WAFOM}(P)$.
\end{enumerate}
In this case, note that the point sets $P$ are not extensible.

\begin{figure}
\centering
  \includegraphics[width=11.5cm]{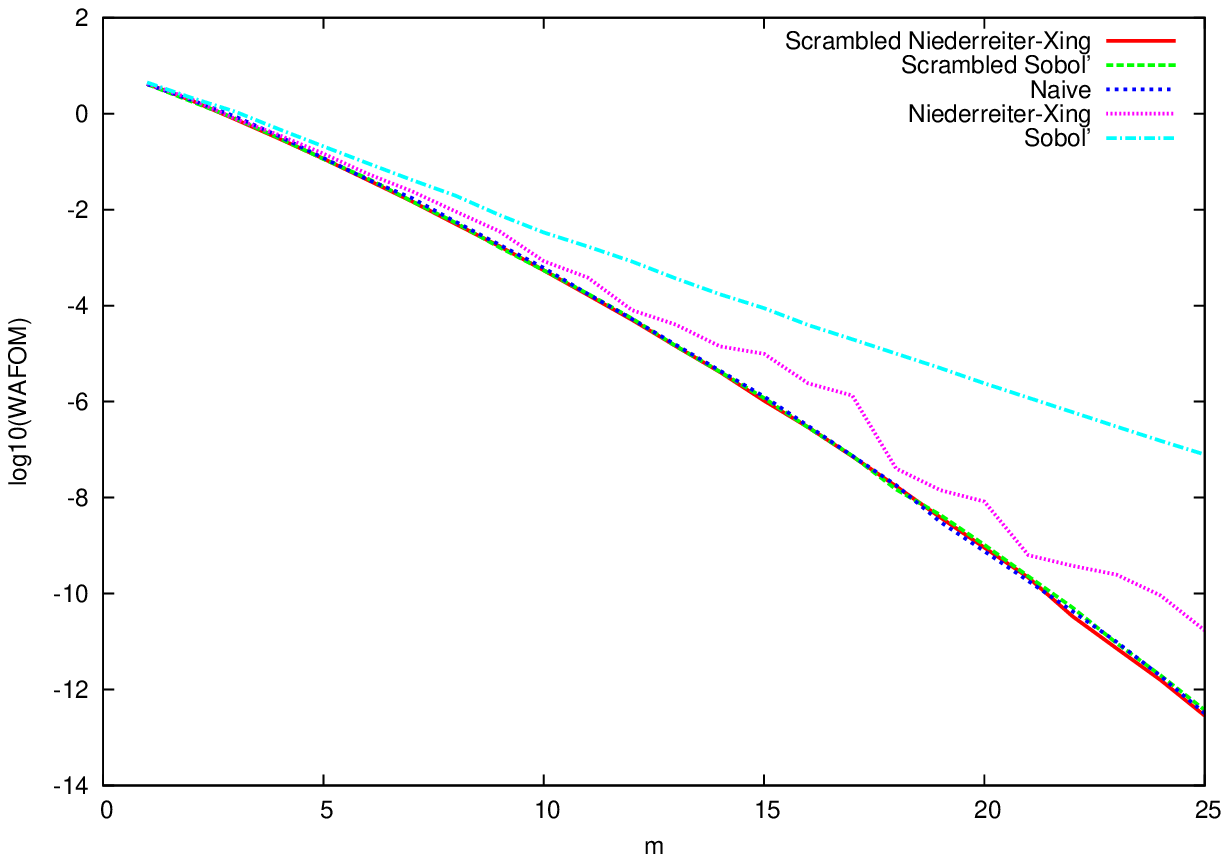}
  \caption{WAFOM (in $\log_{10}$ scale) for $s = 5$ and $m = 1, \ldots, 25$.}
\label{fig:WAFOM}
\end{figure}

As an example, we set  $(s , n, M) = (5, 32, 100000)$ and compare the WAFOM values of the following point sets $P$: 
\begin{enumerate}
\renewcommand{\labelenumi}{(\alph{enumi})}
\item {\bf Niederreiter--Xing} nets \cite{MR1358190} implemented by Pirsic \cite{MR1958872}. 
\item {\bf Sobol'} nets with better two-dimensional projections \cite{MR2429482}.
\item {\bf Naive} low-WAFOM point sets based on a random search \cite{HO}.
\item {\bf Scrambled Niederreiter-Xing} nets given by the above procedure.  
\item {\bf Scrambled Sobol'} nets given by the above procedure. 
\end{enumerate}
Figure~\ref{fig:WAFOM} plots the WAFOM values. 
This shows that (c)--(e) have similar values. 
\textcolor{blue}{The WAFOM values of the Sobol' nets (without linear scrambling) are rather large. 
Roughly speaking, the slope of the Sobol' nets is $O(N^{-1})$. 
Mostly, we can expect the improvement of their efficiency by using linear scrambling. 
Intuitively, we explain these phenomena in terms of WAFOM. 
In (\ref{eqn:WAFOM}), ${\rm WAFOM}(P)$ increases if the proportion of ${x_{i, j}} = 0$ is large. 
(Conversely, ${\rm WAFOM}(P)$ decreases if the proportion of ${x_{i, j}} = 1$ is large.) 
The generating matrices $C_1, \ldots, C_s \in \mathbb{F}_2^{n \times m}$ of the Sobol' nets are 
non-singular upper triangular, and hence the first $2^m$ points always have ${x_{i, j}} = 0$ for $m < j  \leq n$. 
In other words, these least significant bits of the first $2^m$ output points with $n$-digit precision are all zero. 
As a result, ${\rm WAFOM}(P)$ tends to be large in (\ref{eqn:WAFOM}). 
When we apply linear scrambling to the Sobol' nets, 
these least significant bits change from $0$ to $1$ (at random) and the WAFOM values decrease. 
Hence, the rate of convergence is expected to improve. 
On the other hand, the generating matrices of the Niederreiter--Xing nets 
are (almost) dense, and  the WAFOM values are already small, 
so we obtain higher order convergence rates using non-scrambled Niederreiter--Xing nets. 
However, by selecting suitable scrambling matrices, further improvements can be obtained for large values of $m$. 
We conduct additional numerical experiments on these topics in Remark~\ref{remark:high_wafom}. }

\begin{remark} \label{Remark: t-value}
Low-WAFOM point sets based on a simple random search do not always possess small $t$-values, particularly for larger $s$ and $m$. 
Table~\ref{table:t-values} gives a summary of the $t$-values of the above point sets for $s = 5$. 
As described in \cite{HO}, the naive low-WAFOM point sets were
searched by inductively determining the columns vectors of $C_1, \ldots, C_s$ in terms of WAFOM, 
thus allowing extensibility.  
Because we did not consider the \textcolor{blue}{$t$-values} in advance, the $t$-values are rather large. 
Matsumoto--Saito--Matoba (non-extensible) sequential generators \cite{MSM} 
exhibit a similar tendency. 
Nevertheless, such low-WAFOM point sets are effective for smooth functions (see the next section for details). 
\end{remark}
 
\begin{remark}
In two pioneering papers, Dick \cite{MR2346374, MR2391005} proposed {\it higher order digital nets} and {\it sequences} 
that achieve a convergence rate of $O(N^{-\alpha} (\log N )^{\alpha s})$ 
for $\alpha$-smooth functions ($\alpha \geq 1$) by considering the decay of the Walsh coefficients. 
For this, he described an explicit construction for generating matrices, called {\it interlacing}. \textcolor{blue}{First},   
we prepare $s \alpha$ generating matrices $C_1, \ldots, C_{s \alpha} \in \mathbb{F}_2^{m \times m}$ \textcolor{blue}{for} a digital 
$(t, m, s\alpha)$-net in advance. 
These are converted to the matrices $C_1^{(\alpha)}, \ldots, C_s^{(\alpha)} \in \mathbb{F}_2^{{m \alpha} \times m}$ by rearranging the row vectors of  $\alpha$ successive generating matrices. 
Then, the digital net with $C_1^{(\alpha)}, \ldots, C_s^{(\alpha)}$ achieves a convergence rate of $O(N^{-\alpha} (\log N )^{\alpha s})$. 
From \cite[Proposition~15.8]{MR2683394}, such a digital net is a classical digital $(t', m, s)$-net with $t' \leq t$. 
However, when $\alpha$ or $s$ is large, the exact quality parameter $t'$ might become large compared 
with the best possible $t$-value in the framework of classical $(t, m, s)$-nets. 
The last two rows of Table~\ref{table:t-values} give the $t$-values of interlaced 
Niederreiter--Xing nets for $\alpha = 2$ and  $3$. 
Our scrambling approach has the advantages that 
\textcolor{blue}{the exact quality parameter $t$ does} not increase and higher order convergences can be expected. 
\end{remark}

 \begin{table} 

  \scalebox{0.7}{ %
\begin{tabular}{l||c|c|c|c|c|c|c|c|c|c|c|c|c|c|c|c|c|c|c|c|c|c|c|c|c} \hline
  $m$ & $1$ & $2$ & $3$ & $4$ & $5$ & $6$ & $7$ & $8$ & $9$ & $10$ & $11$ & $12$ & $13$ & $14$ & $15$ 
  & $16$ & $17$ & $18$ & $19$ & $20$ & $21$ & $22$ & $23$ & $24$ & $25$ \\ \hline \hline
    {Sobol'} & $0$ & $1$ & $2$ & $2$ & $2$ & $3$ & $3$ & $3$ & $3$ & $3$ & $4$ & $4$ & $5$ & $4$ & $4$ & $5$ & $4$ & $5$ & $5$ & $5$ & $5$ & $5$ & $5$ & $5$ & $5$ \\ \hline 
    {Niederreiter--Xing} & $1$ & $2$ & $1$ & $2$ & $2$ & $2$ & $2$ & $2$ & $2$ & $2$ & $2$ & $2$ & $2$ & $2$ & $2$ & $2$ & $2$ & $2$ & $2$ & $2$ & $2$ & $2$ & $2$ & $2$ & $2$ \\ \hline
    {Naive} & $0$ & $1$ & $2$ & $1$ & $2$ & $3$ & $4$ & $4$ & $4$ & $5$ & $6$ & $7$ & $5$ & $6$ & $6$ & $6$ & $7$ & $7$ & $8$ & $9$ & $9$ & $10$ & $8$ & $9$ & $9$ \\ \hline 
{Interlacing $(\alpha = 2)$} & $1$ & $2$ & $3$ & $4$ & $4$ &  $3$ & $4$ & $4$ & $4$ & $5$ & $6$ & $6$ & $7$ & $6$ & $5$ & $6$ & $7$ & $7$ & $6$ & $6$ & $6$ & $7$ & $6$ & $6$ & $7$ \\ \hline
{Interlacing $(\alpha = 3)$} & $1$ & $2$ & $3$ & $2$ & $3$ & $3$ & $4$ & $5$ & $5$ & $5$ & $6$ & 
$7$ & $6$ & $6$ & $7$ & $8$ & $9$ & $9$ & $7$ & $8$ & $8$ & $9$ & $8$ & $8$ & $8$ \\ \hline
\end{tabular}
}
 \caption{The exact quality parameters $t$ for $m = 1, \ldots, 25$ and $s = 5$.}
\label{table:t-values}
\end{table}

\begin{remark}
Goda, Ohori, Suzuki, and Yoshiki \cite{GOSY} proposed a variant of WAFOM from the \textcolor{blue}{viewpoint} of the mean square error for digitally shifted digital nets. 
They defined the criterion by 
replacing $2$ in (\ref{eqn:WAFOM}) with $4$. Thus, this is similarly applicable to our approach.
\end{remark}

\section{Numerical results} \label{Sec:integration}

To evaluate the point sets (a)--(e) described in Section~\ref{Sec:scrambling}, 
we applied the Genz test package \cite{Genz1984,Genz1987}.
This has been used in many studies (e.g., \cite{MR1417864, MR1958872,mSLO94a, MR1849865}), and was also analyzed from a theoretical
perspective in \cite{MR1963917}. 
Thus, we investigate six different test functions defined over $[0,1)^s$. These are:
\begin{eqnarray*}
\begin{array}{ll}
\mbox{Oscillatory:} & f_1(\mathbf{x})  =  \cos (2 \pi u_1 + \sum_{i = 1}^{s} a_i x_i),\\
\mbox{Product Peak:} & f_2(\mathbf{x})  =  \prod_{i = 1}^{s} [1/{(a_i^{-2} + (x_i - u_i)^2)]},\\
\mbox{Corner Peak:} & f_3(\mathbf{x})  =  (1 + \sum_{i = 1}^{s} a_i x_i)^{-(s+1)},\\
\mbox{Gaussian:} & f_4(\mathbf{x})  =  \exp (-\sum_{i=1}^s a_i^2 (x_i - u_i)^2),\\
\mbox{Continuous:} & f_5(\mathbf{x})  =  \exp (-\sum_{i=1}^s a_i |x_i - u_i|),\\
\mbox{Discontinuous:} & f_6(\mathbf{x})  =  \left\{
\begin{array}{ll}
0, & \mbox{if $x_1 > u_1$ or $x_2 > u_2$},\\
\exp (\sum_{i = 1}^s a_i x_i), & \mbox{otherwise.}
\end{array}
\right.
\end{array}
\end{eqnarray*}
In these functions, we have two parameters, i.e., the difficulty parameters $\mathbf{a} = (a_1, \ldots, a_s)$ and the shift parameters $\mathbf{u} = (u_1, \ldots, u_s)$. 
We generate $\mathbf{a}= (a_1, \ldots, a_s) $ and $\mathbf{u}= (u_1, \ldots, u_s)$ as uniform random vectors in $[0, 1]^s$, and renormalized $\mathbf{a}$ to satisfy the following condition:
\begin{eqnarray*} \label{eqn:Genz condition}
\sum_{i = 1}^s a_i = h_j,
\end{eqnarray*}
where $h_j$ depends on the family $f_j$. 
By varying $\mathbf{a}$ and $\mathbf{u}$,  
we formed quantitative examples based on 20 random samples for each function class. 
For any sample size $|P|=2^m$ and any function $f_j$, 
we computed the median of the relative errors (in $\log_{10}$ scale)
\[ \log_{10} \frac{|I(f_j) - I_N(f_j)|}{|I(f_j)|} \]
varying the parameters, where $I(f_j) := \int_{[0,1)^s} f_j \textcolor{blue}{(\mathbf{x})} \textrm{d} \mathbf{x}$, $N := |P|$, and $I_N(f_j) := ({1}/{|P|}) \sum_{\mathbf{x} \in P} \textcolor{blue}{f_j(\mathbf{x})}$. 

\begin{figure}
  \includegraphics[width=8cm]{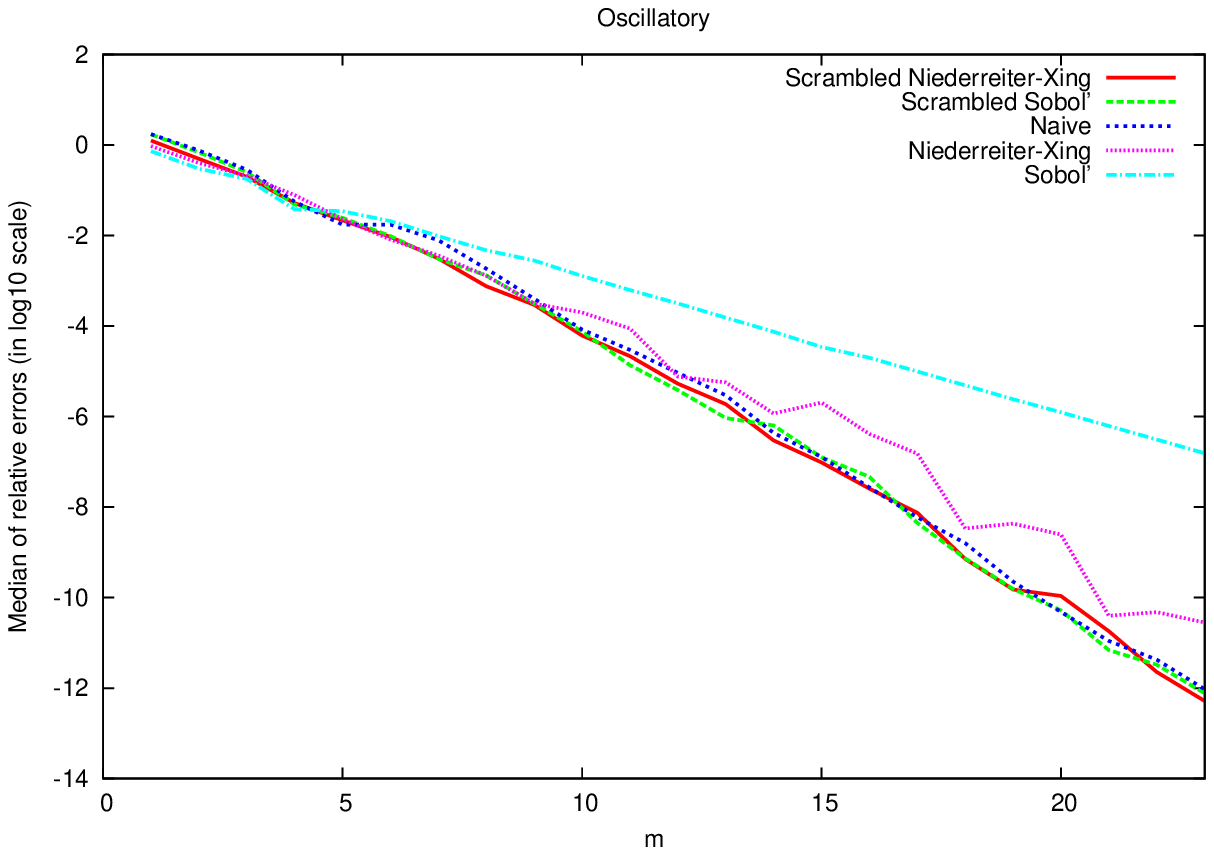}
  \includegraphics[width=8cm]{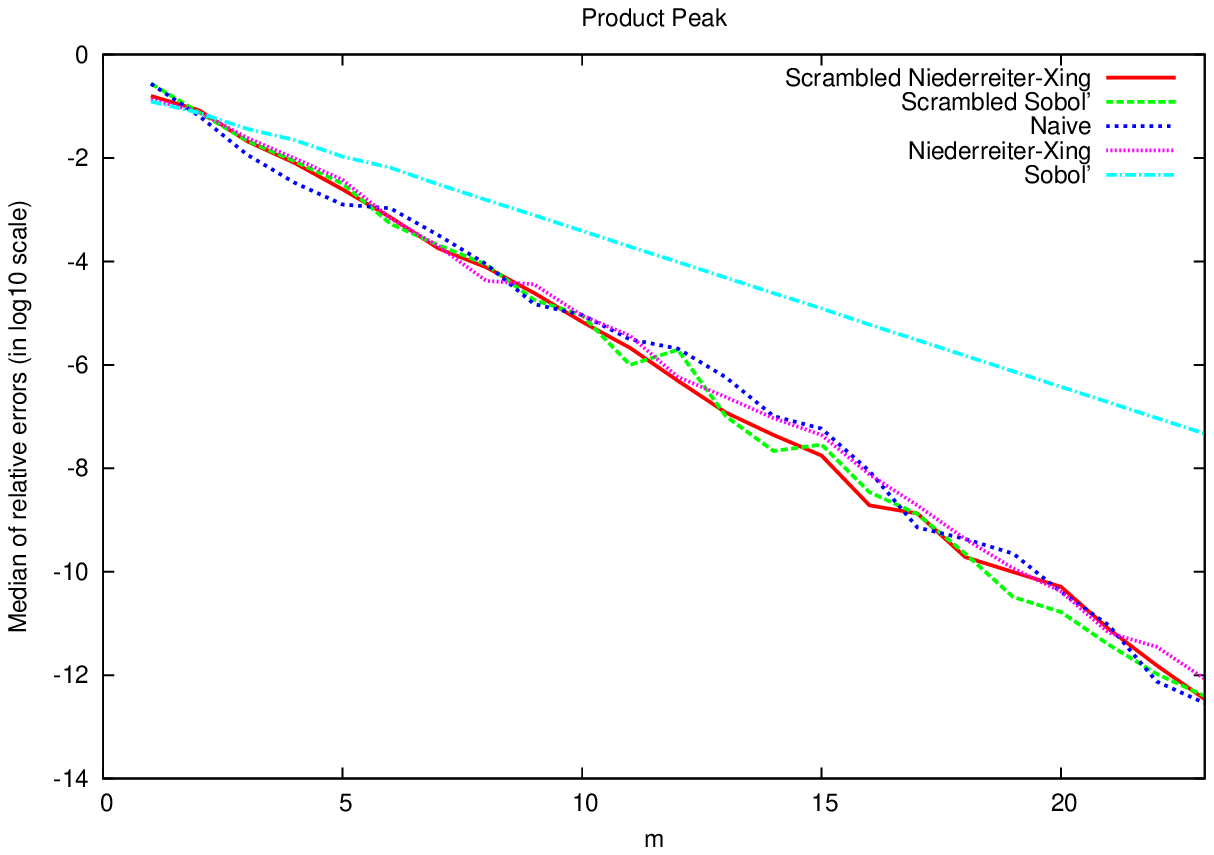}
  \includegraphics[width=8cm]{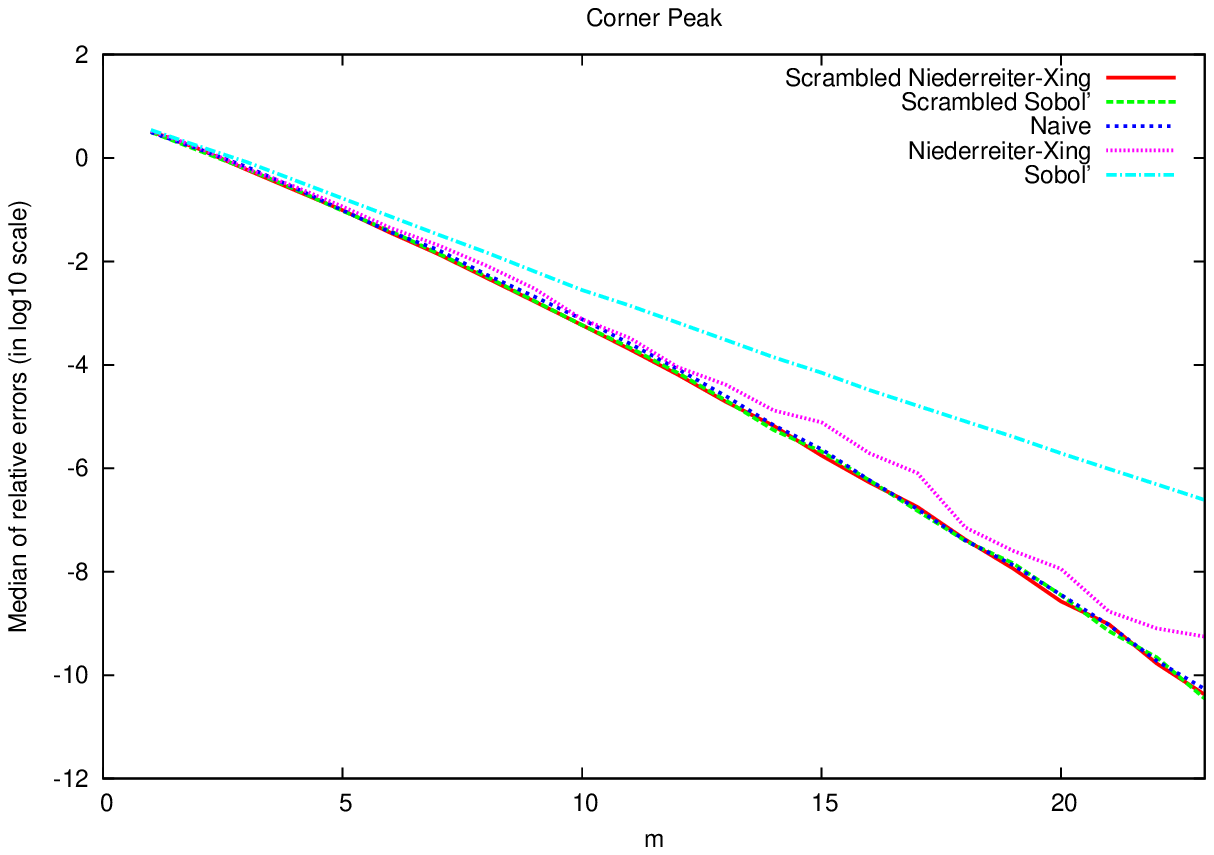}
  \includegraphics[width=8cm]{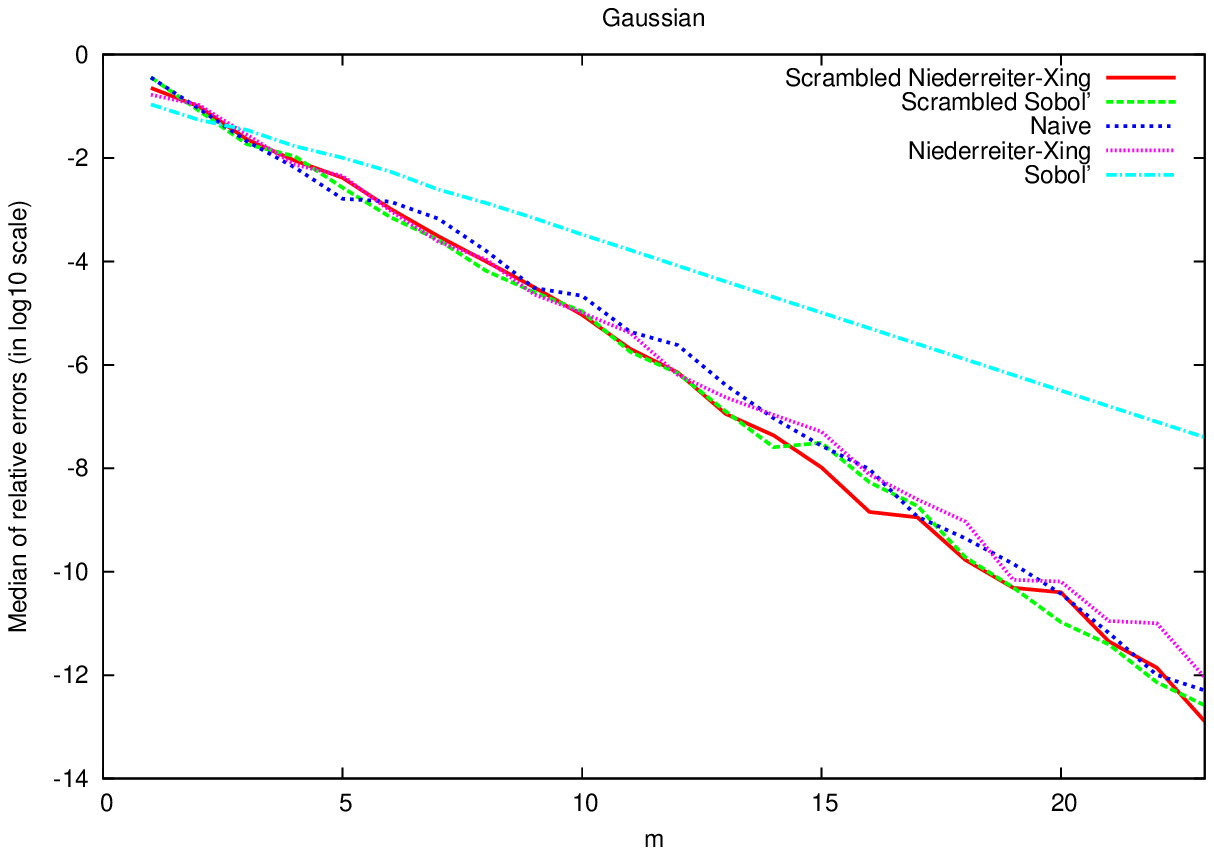}
  \includegraphics[width=8cm]{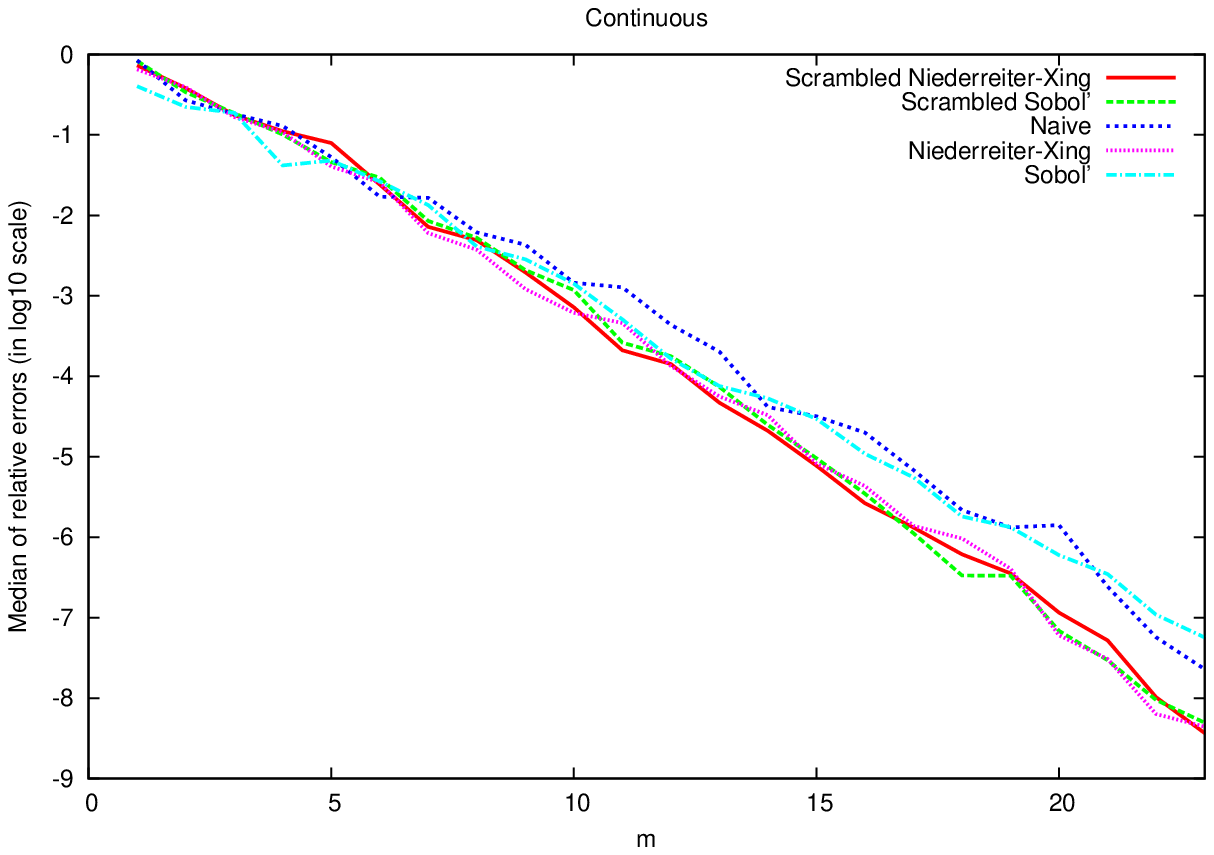}
  \includegraphics[width=8cm]{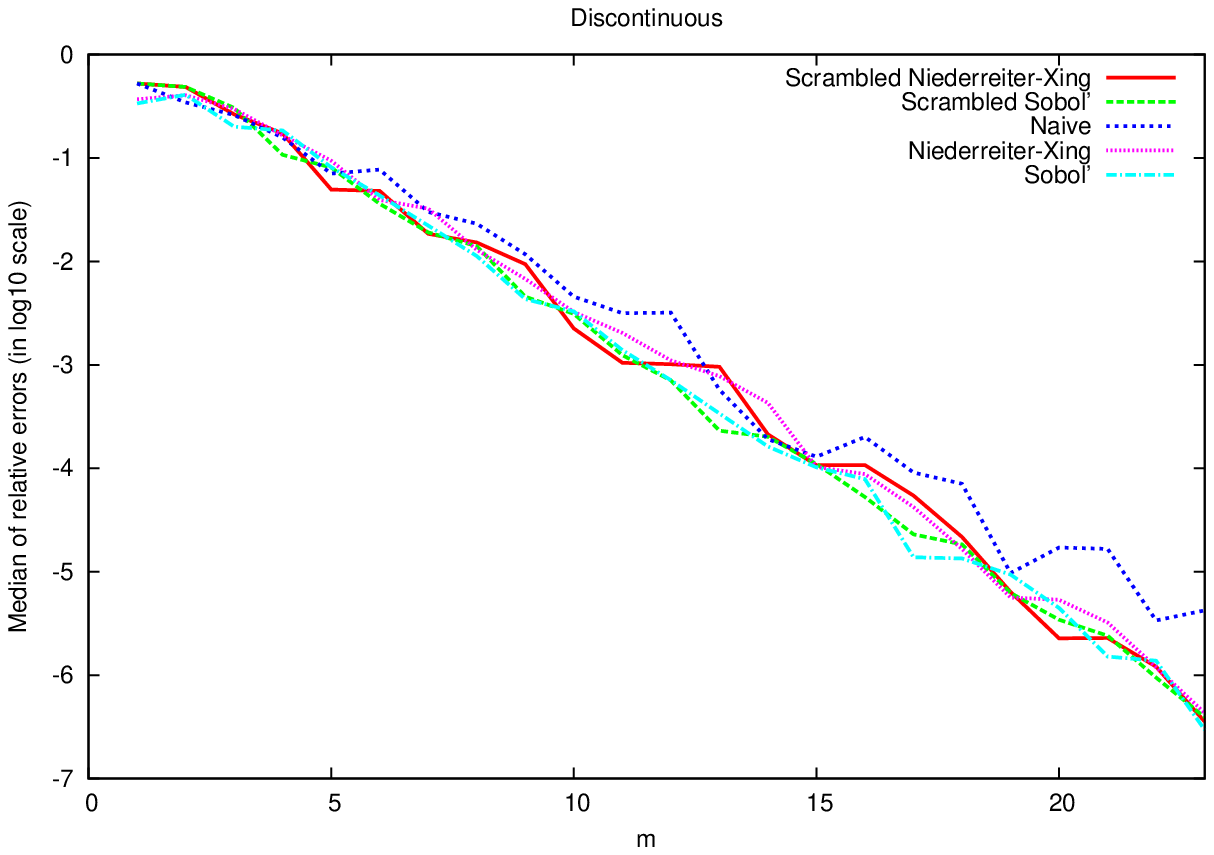}
    \caption{Median of relative errors for Genz functions.}
    \label{fig:Genz}
\end{figure}

Figure~\ref{fig:Genz} shows a summary of the medians of the relative errors 
for $s = 5$, $m = 1, \ldots, 23$, and $(h_1, \ldots, h_6) = (4.5, 3.625, 0.925, 3.515, \textcolor{blue}{10.2}, 2.15)$, 
which are \textcolor{blue}{the settings used} in \cite{HO}. 
For $f_1$and $f_3$, the low-WAFOM point sets 
are clearly superior to the Niederreiter-Xing nets. 
In particular, the scrambled Sobol' nets represent a drastic improvement over the original Sobol' nets. 
Note that the slopes are similar to those in Figure~\ref{fig:WAFOM}. 
Additionally, for $f_2$ and $f_4$, the low-WAFOM point sets are competitive with the Niederreiter--Xing nets. 
In these smooth functions, the WAFOM criterion seems to work very well. 
In the case of non-smooth functions, the situations \textcolor{blue}{are} different. 
For the continuous but non-differentiable functions $f_5$,
the naive low-WAFOM point sets are inferior to the Niederreiter--Xing nets. 
However, when we take into account the \textcolor{blue}{$t$-value} of $(t, m, s)$-nets, 
the low-WAFOM point sets preserve the rate of convergence. 
For $f_6$, 
the naive low-WAFOM point sets are also inferior to the other point sets with small $t$-values. 
These results imply that the \textcolor{blue}{$t$-value} is important for non-smooth functions. 

\textcolor{blue}{Finally, we note that, as the dimension $s$ increases, the WAFOM values tend to have only slight differences (see Section~4.2 of \cite{MO2014} for details). 
In this case, the rates of convergence weaken, but the obtained point sets in this paper seem to be at worst comparable to the original non-scrambled Niederreiter--Xing or Sobol' nets, especially for high-smooth functions. (To save space, we omit the figures.)}

\begin{remark} \label{remark:high_wafom} 
\textcolor{blue}{
There are some experimental reports that random linear scrambling improves the rates of convergence in numerical integration. 
To investigate the effect of  WAFOM and scrambling,  
we conduct further experiments on a comparison between 
scrambled nets with small WAFOM and those with large WAFOM. 
For this purpose, using the similar algorithm to that in Section~\ref{Sec:scrambling}, we searched for linearly scrambled digital $(t, m ,s)$-nets $P$ with small $t$-values but
with the largest ${\rm WAFOM}(P)$:
\begin{enumerate}
\renewcommand{\labelenumi}{(\alph{enumi})}
\setcounter{enumi}{5}
\item {\bf Scrambled Niederreiter--Xing (worst)} nets with the largest ${\rm WAFOM}(P)$.  
\item {\bf Scrambled Sobol' (worst)} nets with the largest ${\rm WAFOM}(P)$. 
\end{enumerate}
Figure~\ref{fig:Genz2} plots the WAFOM values and the medians of relative errors of the Genz function packages for
the point sets (a), (b), and (d)--(g) in the same settings as in Figure~\ref{fig:Genz}. 
{\bf Scrambled Niederreiter--Xing (best)} and {\bf Scrambled Sobol' (best)} are copies of (d) and (e) in Figure~\ref{fig:Genz} (with the smallest ${\rm WAFOM}(P)$), respectively.}
\textcolor{blue}{
We can summarize our experimental results as follows:
\begin{itemize}
\item The largest WAFOM values of the scrambled Sobol' nets are 
comparable to or slightly better than the WAFOM values of the non-scrambled Sobol' nets. 
Thus, most scrambled Sobol' nets have WAFOM values that are smaller than those of the non-scrambled Sobol' nets (as pointed out in Section~\ref{Sec:scrambling}), 
and hence we can expect that the simple application of ``random" linear scrambling improves the rate of convergence for the Sobol' nets from the viewpoint of WAFOM. 
In Figure~\ref{fig:Genz2}, 
the scrambled Sobol' nets with the largest WAFOM are better than the non-scrambled Sobol' nets for all the smooth functions, especially $f_2$ and $f_4$,  
but the scrambled Sobol' nets with the smallest WAFOM seem to be the best choices. 
\item The WAFOM values of the Niederreiter--Xing nets are already small, 
and the WAFOM values of the scrambled Niederreiter--Xing nets given by inappropriate lower triangular matrices become larger than 
those of the non-scrambled Niederreiter--Xing nets. 
Indeed, 
the scrambled Niederreiter--Xing nets with the largest WAFOM are worse than the non-scrambled Niederreiter--Xing nets for all the smooth Genz functions. 
\end{itemize}
Overall, WAFOM is a good criterion for ensuring higher order convergence for high-smooth functions.
}


\begin{figure}
\begin{center}
\includegraphics[width=7cm]{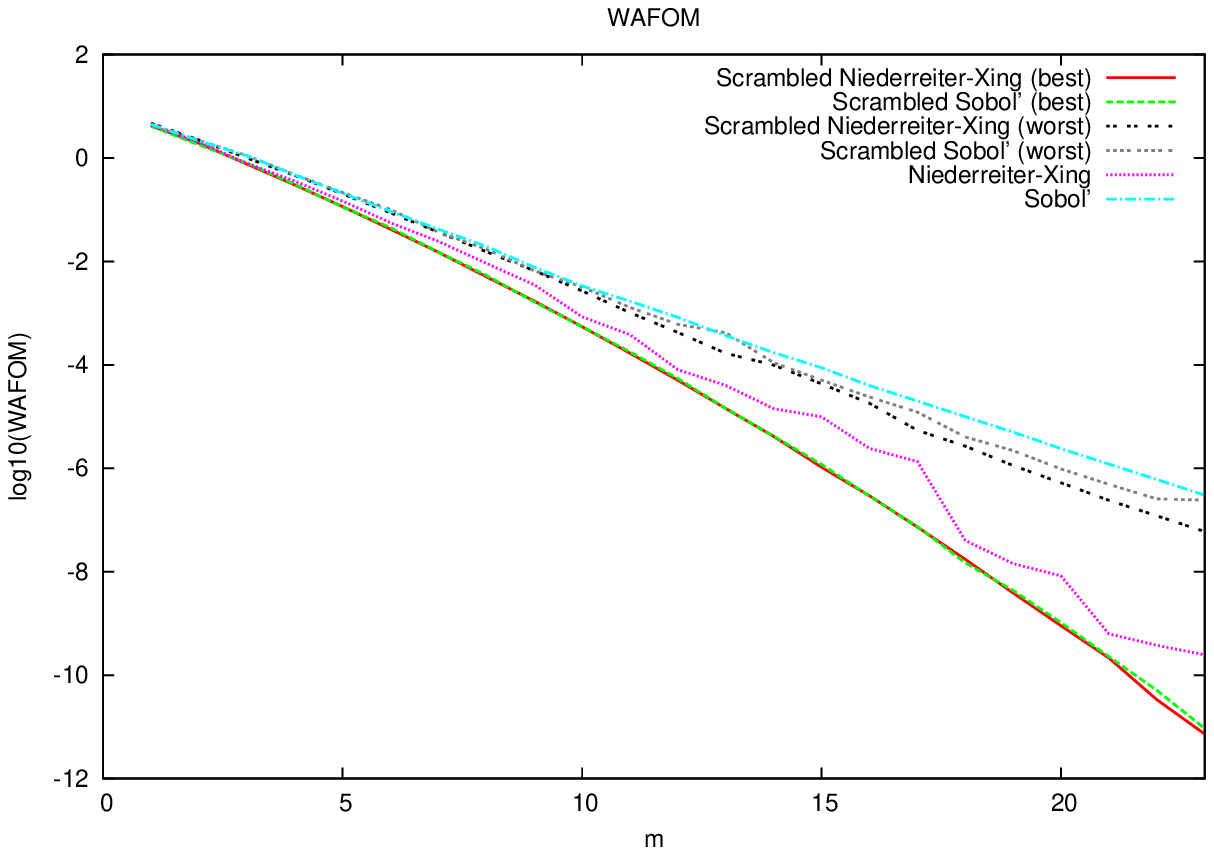}\\
\end{center}
\includegraphics[width=7cm]{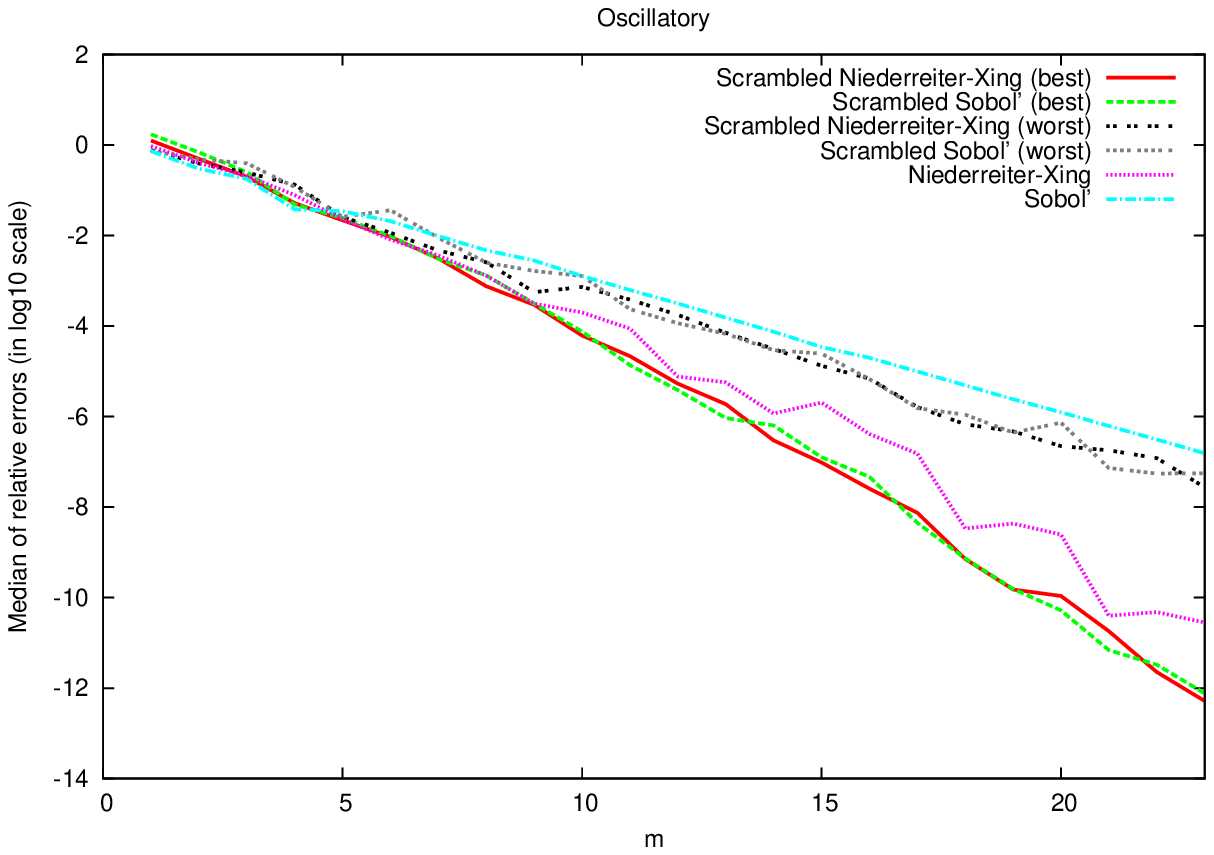}
\includegraphics[width=7cm]{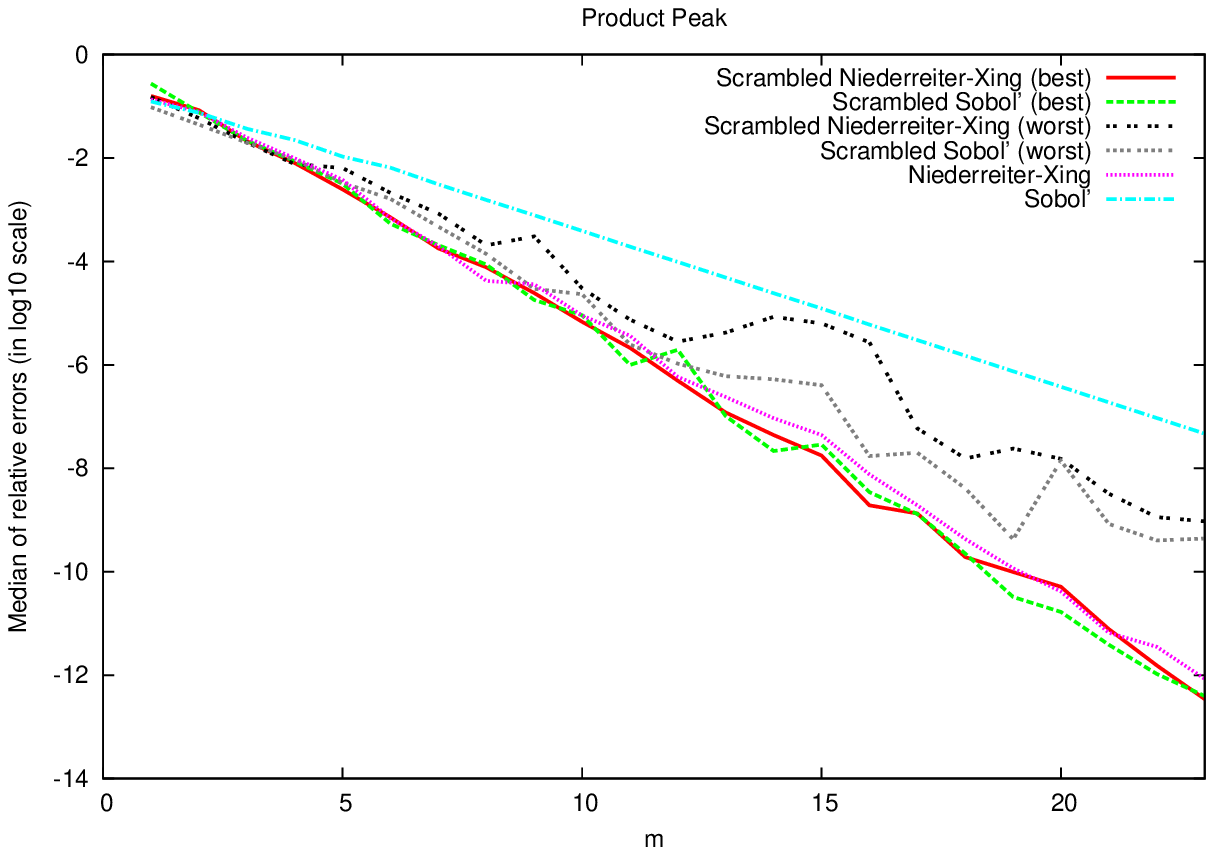}
\includegraphics[width=7cm]{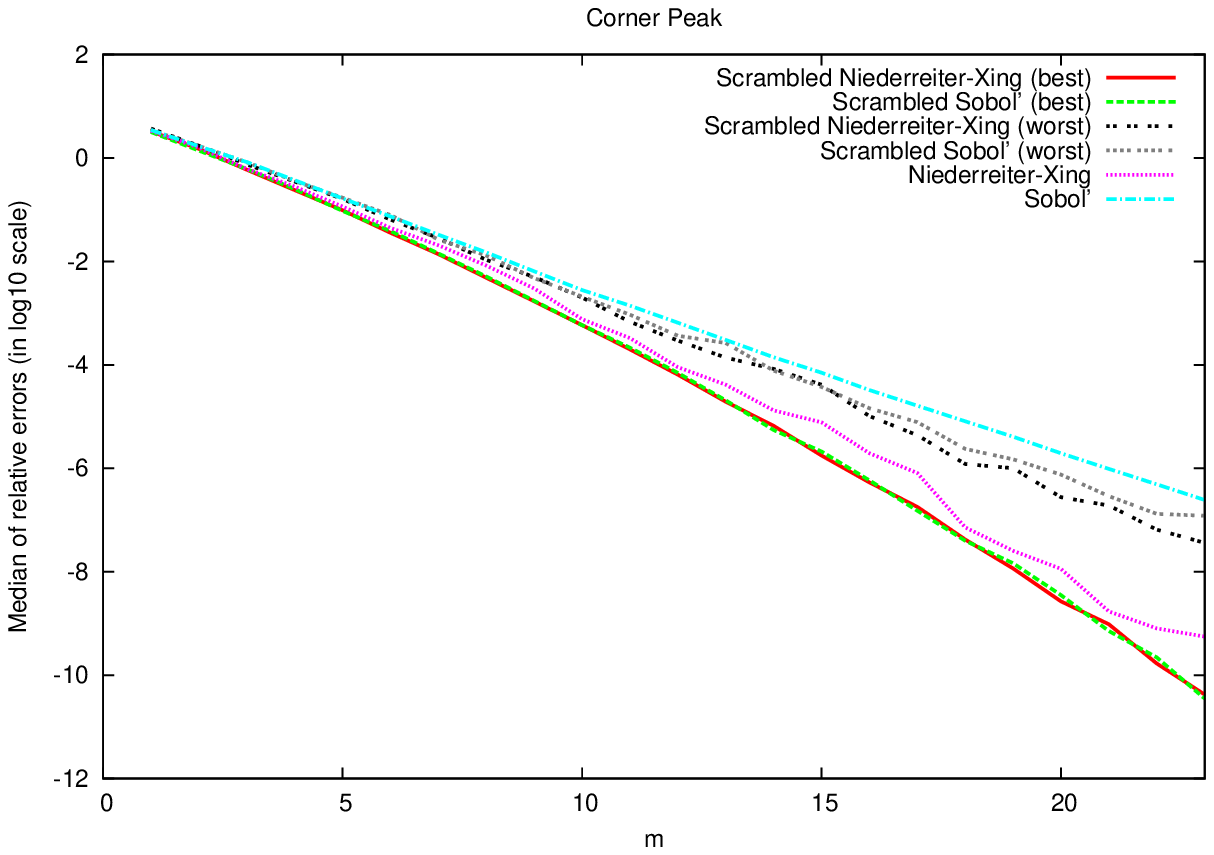}
\includegraphics[width=7cm]{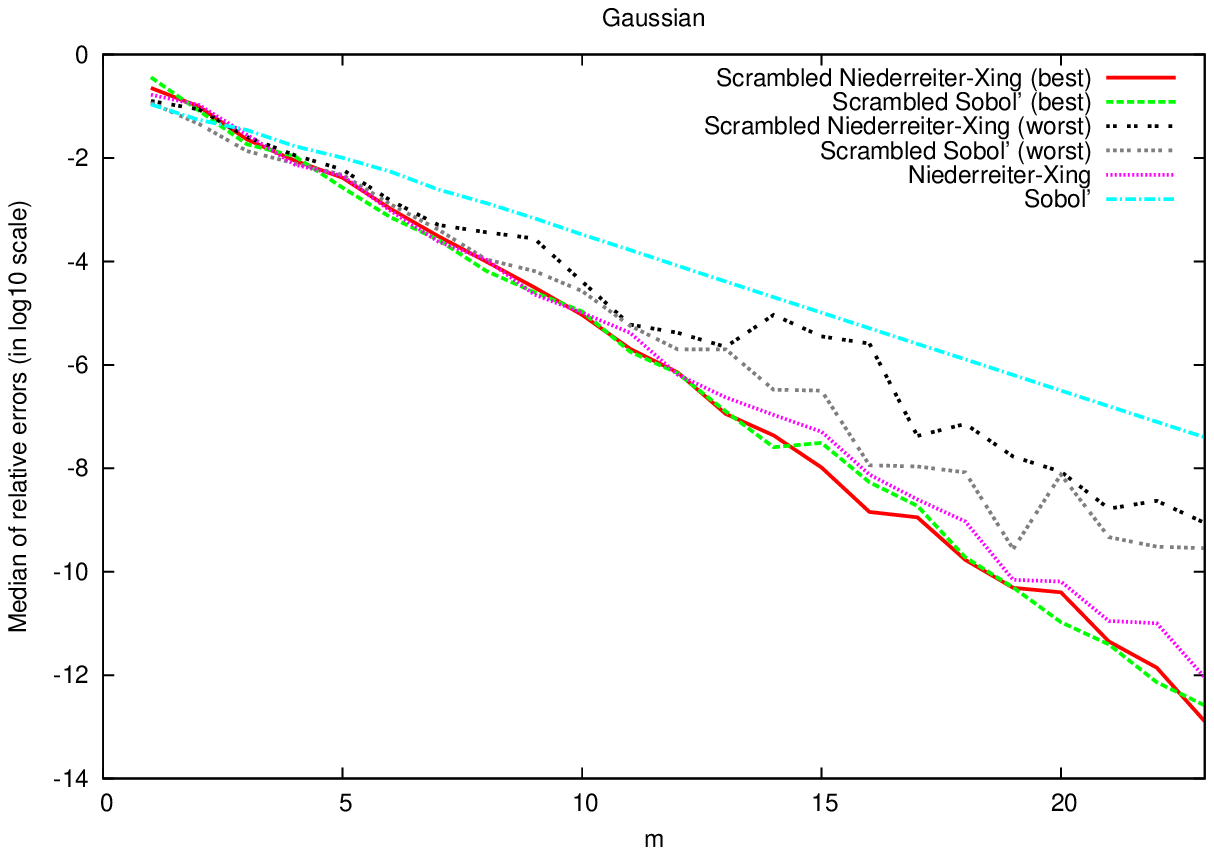}
\includegraphics[width=7cm]{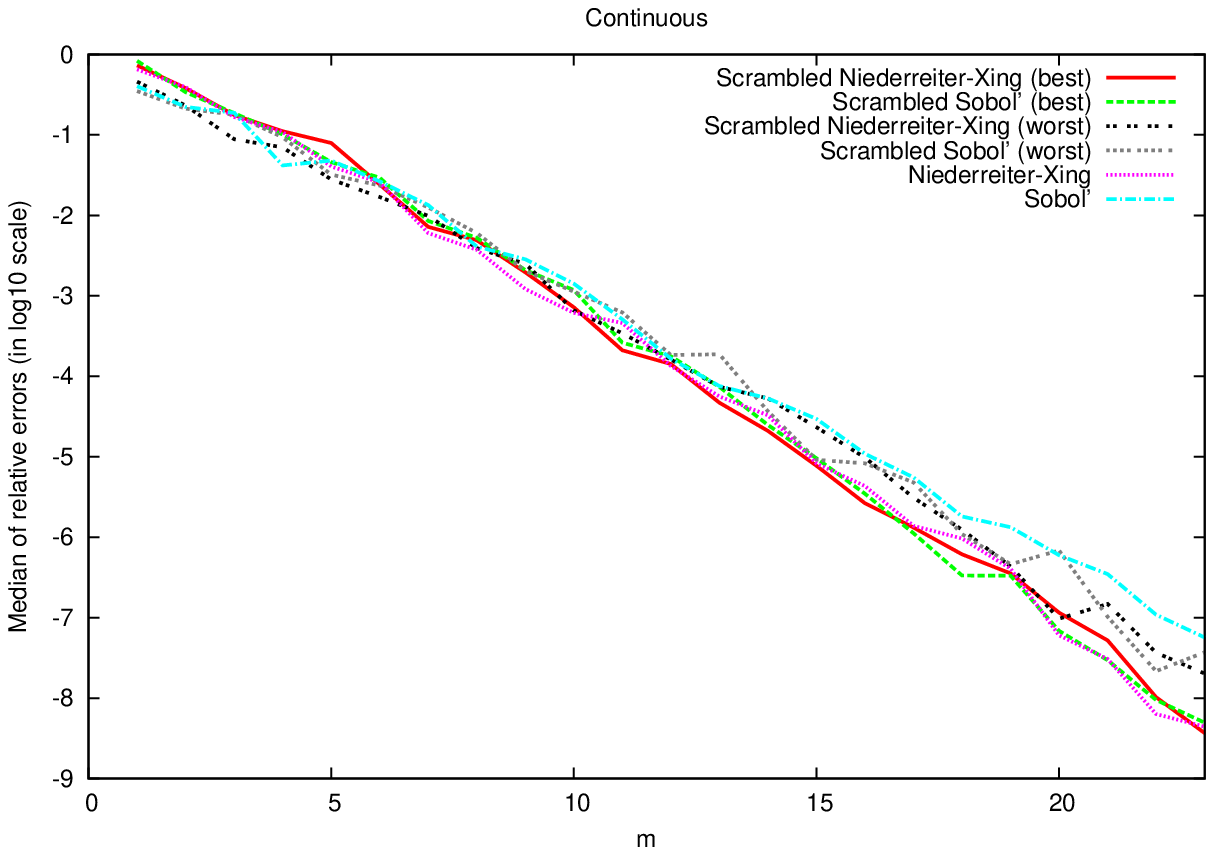}
\includegraphics[width=7cm]{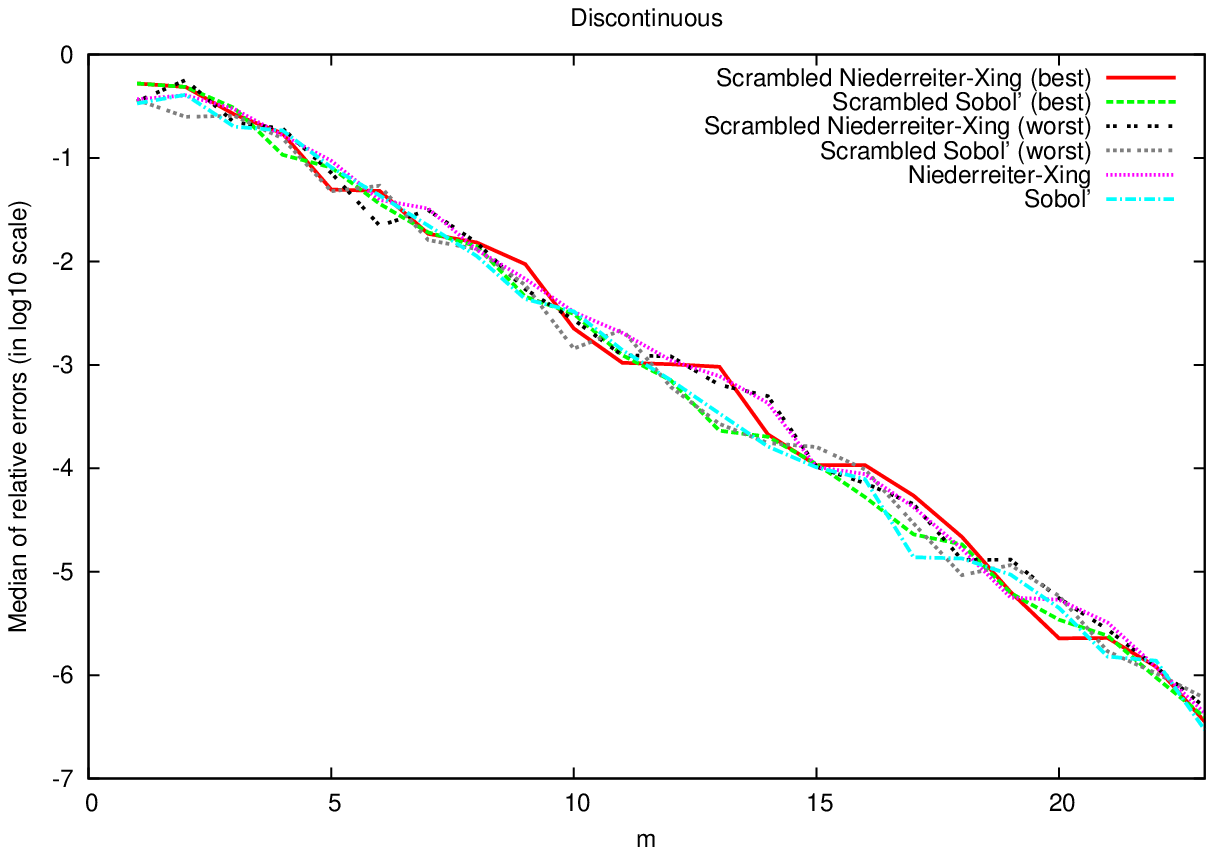}
\caption{\textcolor{blue}{Comparison of scrambled digital nets with small WAFOM and those with large WAFOM for $s = 5$. The top figure shows WAFOM values (in $\log_{10}$ scale) for $m = 1, \ldots, 23$. 
The other figures show the median of relative errors for the Genz functions for $m = 1, \ldots, 23$. }}
\label{fig:Genz2}
\end{figure} 
\end{remark}

\section{Conclusions and future directions}

In this paper, we have searched for  point sets whose $t$-value and WAFOM are both small
so as to be effective for a wider range of function classes. 
For this, we fixed digital $(t, m, s)$-nets in advance and applied random linear scrambling. 
The key technique was \textcolor{blue}{the selection} of linearly scrambled $(t, m, s)$-nets in terms of WAFOM. 
Numerical experiments showed that the point sets obtained by our method have improved convergence rates for smooth functions and are robust for non-smooth functions. 

Finally, we discuss some directions for future research.
In our approach, $m$ was fixed and the extensibility was discarded. 
We also attempted to search for extensible point sets, 
but the WAFOM values tended to be worse than the current ones for large $m$. 
Thus, an efficient search algorithm for extensible scrambling matrices is 
one area of future work. 
As another direction, the quasi-Monte Carlo method is an important tool in computational finance (e.g., \cite{MR1999614, MR2519835}). 
However, many applications encounter integrands with boundary singularities. 
Such integrands are not included in a suitable class of functions, i.e., $n$-smooth functions,  
so we might not expect higher order convergence from the simple application of low-WAFOM point sets. 
There will probably be a need for some kind of transformation to force the integrand to be included in a suitable class of functions, 
such as periodization in lattice rules. 
The study of WAFOM is still in its infancy, 
so a number of unsolved problems remain. 

\subsection*{Acknowledgments}
\textcolor{blue}{
The author is thankful to the
anonymous referees for their valuable comments and suggestions.
The author also wishes to express his gratitude to Professor Makoto Matsumoto at Hiroshima University and Professor Syoiti Ninomiya at Tokyo Institute of Technology for continuous encouragement and many helpful comments.
The author was partially supported by Grant-in-Aid for JSPS Fellows 24$\cdot$7985, Young Scientists (B) 80610576, and Scientific Research (B) 70231602.}





\bibliographystyle{model1b-num-names}
\bibliography{wafombib}







\end{document}